\documentclass[12pt]{article}
\usepackage{amsmath,amssymb}
\usepackage[dvips]{graphicx}
\newtheorem{lem}{Lemma}[section]
\newtheorem{thm}[lem]{Theorem}
\newtheorem{prop}[lem]{Proposition}

\newtheorem{Def}[lem]{Definition}
\newcommand{\beqn}{\begin{eqnarray}}
\newcommand{\eeqn}{\end{eqnarray}}

\newcommand{\ep}{\varepsilon}

\newcommand{\e}{\pmb{e}}

\newcommand{\C}{{\mathbb{C}}}
\newcommand{\Z}{{\mathbb{Z}}}

\newcommand{\Pf}{\mathrm{Pf}}

\newcommand{\diag}{{\mathrm{diag}}}

\newcommand{\la}{\lambda}

\begin{document}
\title
{{\sf 
Schubert classes in the equivariant
cohomology of \\
the Lagrangian Grassmannian
}
}
\date{} 
\author
{{\scshape Takeshi Ikeda }
}

\bigskip

\maketitle

\begin{center}
{\it
To Ryoshi Hotta on the occasion of his 65th birthday}
\end{center}

\begin{abstract}
Let $LG_n$ denote the Lagrangian Grassmannian
parametrizing maximal isotropic (Lagrangian)
subspaces of a fixed symplectic vector space of dimension $2n.$
For each strict partition $\lambda=(\lambda_1,\ldots,\lambda_k)$ with 
$\lambda_1\leq n$ there is a Schubert variety $X(\lambda).$
Let $T$ denote a maximal torus of the symplectic group
acting on $LG_n.$
Consider the $T$-equivariant cohomology 
of $LG_n$ and the
$T$-equivariant fundamental class $\sigma(\lambda)$
of $X(\lambda).$
The main result of the present paper is an explicit formula
for the restriction of the
class $\sigma(\lambda)$
to any torus fixed point.
The formula is written in terms of 
factorial analogue of the Schur $Q$-function,
introduced by Ivanov.
As a corollary to the restriction formula,
we obtain an equivariant version
of the Giambelli-type formula for $LG_n.$
As another consequence of the main result, 
we obtained a presentation of the
ring $H_T^*(LG_n).$ 
\end{abstract}

\section{Introduction}
\setcounter{equation}{0}

\bigskip

Let $LG_n$ denote 
the \textit{Lagrangian Grassmannian\/}
parametrizing
$n$-dimensional
isotropic subspaces of a fixed $2n$-dimensional
symplectic vector space.
The \textit{Schubert classes} give a linear basis
for the integral cohomology ring of $LG_n.$
These classes can be parametrized by 
sequences $\la=(\la_1,\ldots,\la_k)$ of 
integers
such that $n\geq \la_1>\cdots>\la_k\geq 1.$ 
The same set of sequences also 
parametrizes the $T$-fixed points in $LG_n,$
where $T$ denotes a maximal torus of the symplectic
group acting on $LG_n.$
Here we consider the 
$T$-equivariant integral cohomology ring $H_T^*(LG_n).$
We are interested in the 
$T$-\textit{equivariant 
Schubert classes} $\sigma(\la)$ in $H_T^*(LG_n).$ 
The classes $\sigma(\la)$ form a free basis 
of $H_T^*(LG_n)$ over the ring $\mathcal{S}$ 
of $T$-equivariant cohomology of a point.
It is known that $\mathcal{S}$ is naturally identified with the polynomial ring
$\Z[\ep_1,\ldots,\ep_n]$
(for the definition of $\ep_i$ see \S \ref{ESC}).
Let $e(\mu)$ be the $T$-fixed point corresponding to 
an index $\mu.$
The inclusion map $i_\mu:\{e(\mu)\}\hookrightarrow LG_n$ 
induces the \textit{restriction} morphism
$i_\mu^*:H_T^*(LG_n)\rightarrow H_T^*(\{e(\mu)\})\cong \mathcal{S}.$
The main result of the present paper (Theorem \ref{res}) is
an explicit formula  
for $i_\mu^*\sigma(\la)$, the restriction of the equivariant 
Schubert class $\sigma(\la)$ to a $T$-fixed point $e(\mu),$
as a polynomial in $\Z[\ep_1,\ldots,\ep_n].$ 

For the classical 
Grassmannian $G_{d,n}$ of 
$d$-dimensional subspaces of $\C^n$,
Knutson and Tao \cite{KT} derived
a formula for the restriction of the equivarent Schubert
classes to a torus fixed point. 
This formula is written in terms of 
the \textit{factorial analogue} of 
the Schur function introduced by Biedenharn and Louck \cite{BL}
and studied by other authors (see Section \ref{Fac}). 
Recently, Lakshmibai, Raghavan, and Sankaran \cite{LRS}
also proved the restriction formula 
by a different method. 
They also derived 
\textit{equivariant Giambelli formulas},
determinantal formulas in the global 
ring $H_T^*(G_{d,n})$ of equivariant cohomology,
that express the 
equivariant Schubert classes.

As for the ordinary cohomology
of the Lagrangian Grassmannian, 
Hiller and Boe \cite{HB} proved 
Pieri-type formulas.
Using the result in \cite{HB}, Pragacz \cite{P} proved Giambelli-type formulas
that express each Schubert class in a Schur-type Pfaffian form.
The derivation of the Pfaffian formula 
is based on a comparison of the formula in \cite{HB}
and the Pieri formula for the Schur $Q$-function.
Our formula (Theorem \ref{res}) is written in terms of 
the factorial analogue 
of the Schur $Q$-function, introduced by Ivanov \cite{Iv1,Iv2}.
This 
leads to an 
equivariant Giambelli-type formula (\ref{Giambelli})  
analogous to Pragacz's result.

Here is a brief summary of the paper. In Section \ref{Pre} we
recall some results on the Weyl group of type $C_n$
and fix some standard notation.
The set of minimal coset representives 
introduced in this section will be used as an index set
that parametrize the main objects of this paper,
the Schubert varieties, $T$-fixed points etc.
We also present several combinatorial 
descriptions of this set. 
Section \ref{Lag} is devoted to 
basic geometric settings,
where we introduce the Schubert vatieties.
We proceed by studying the $T$-equivariant cohomology 
of $LG_n$ in Section \ref{Equ}.
A recurrence relation (\ref{Rec}) arising from 
the equivariant Pieri-Chevalley formula (\ref{Ch})
plays the central role.
In Section \ref{Fac} we give a definition of the 
factorial analogue of the Schur $Q$-function and 
present some properties. 
Finally, in Section \ref{Main} we prove
our main theorem.
The proof is performed by a comparison 
of Ivanov's Pieri-type formula (\ref{Pieri}) for factorial
$Q$-function and the recurrence relation for the 
restricted Schubert class.
Using the restriction formula, 
we prove the Giambelli-type formula (\ref{Giambelli})
for the equivariant Schubert classes of 
$LG_n.$
In Section \ref{Two},
we also give a supplementary 
discussion on expressing the Schubert 
classes of {\it two-row\/} type diagrams in terms 
of a polynomial in the {\it special\/} Schubert
classes $\sigma(i)\;(1\leq i\leq n).$
A formula in Prop. \ref{expansion} seems to be new. It
provides an explicit form of any two-row type 
factorial $Q$-function.
In Section \ref{Ring}, we prove 
a ring presentation for $H_T^*(LG_n)$ as 
a quotient of the polynomial ring over $\mathcal{S}.$
In Appendix (Section \ref{App}), we present a 
brief introduction to Ivanov's functions 
and prove a vanishing property crucial 
to the main body of the paper.


It is well-known that the cohomology ring
of the Grassmannian of orthogonal type
is very similar to that of $LG_n.$
It is natural to expect
the equivariant Schubert classes for the 
orthogonal Grassmannian
is described by factorial Schur $P$-functions (cf. \cite{P},\cite{Jo}).
Indeed, after the first version of this paper
is written, the author and H. Naruse succeed in 
deriving such formulas. We will discuss the 
subject in a separate paper \cite{IN}
which focus on some  
combinatorial aspects of
the equivariant Schubert calculus.

It should be mentioned that 
a recent result due to 
Ghorpade and Raghavan \cite{GR}
provides an alternative combinatorial approach to our formula (\ref{MAIN}).
They presented the coordinate ring of the tangent cone of $X(\lambda)$ 
at $e(\mu)$ as a Stanley-Reisner ring.
Note that a preceding result of Conca \cite{Co}
corresponds to the case with the fixed point
is the identity coset.
This description immediately leads to a combinatorial
formula for $\sigma(\lambda)|_\mu$,
which 
is quite similar to a tableau type formula for $Q_\lambda(x|a)$
given in \cite{Iv2}. Details of these issues will be 
discussed in \cite{IN}.
Recently, the result, namely a combinatorial
expression for the restriction
of a Schubert class in $LG_n$
to a $T$-fixed point,
was independently proved by Kreiman \cite{Kr2}.
%
\bigskip

\textit{Acknowledgments}: Firstly, I would like to thank 
H. Naruse for his keen interest
in the present work and plenty of helpful comments.
I also thank K. Takasaki, Y. Kodama, T. Tanisaki,
and M. Kaneda for valuable discussions.
During the prepation of this paper 
I also benefited from coversations
with M. Ishikawa, S. Kakei, M. Katori, and T. Maeno. 
I also thank H. Nakajima and the referee who encouraged me  
to prove a Giambelli-type formula (\ref{Giambelli}) which was
not included in the first version of this paper.
Lastly, but not leastly, I am grateful to
H.-F. Yamada for showing me the importance of Schur's $Q$-functions.
This research was partially supported by 
Grant-in-Aids for
Young Scientists (B) (No. 17740101)
from Japan Society of
the Promotion of Science.

\section{Preliminaries}\label{Pre}
\setcounter{equation}{0}
We first recall some basic notions about 
the Weyl group of type $C_n$, in order to fix our notation.
The purpose of this section is to introduce 
an index set for the main
ingredients of this paper, Schubert varieties, torus fixed points, etc.
References for this section 
are \cite{B},\cite{H}, and \cite{HB}.
	\subsection{Weyl group of type $C_n$ 
	}
Let $S_{2n}$ be the symmetric group
of all permutations of $2n$ 
letters $\{1,\ldots,2n\}.$
Set $\overline{i}=2n-i+1.$
Let $W$ be the subgroup of $w\in S_{2n}$ such that 
$w(i)=j\Longleftrightarrow w(\bar{i})=\bar{j}.$
Then $w\in W$ can be determined by $w(1),\ldots,w(n).$
A standard set of generators of $W$ is given by 
$s_i=(i,i+1)(\overline{i+1},\overline{i})\;
(1\leq i\leq n-1)$ and $s_n=(n,\overline{n}),$
where we denote by $(i,j)$ the transposition.
The \textit{length} $\ell(w)$ of an element $w$ in
$W$ is the smallest number of the generators $s_1,\ldots,s_n$
(the simple reflections)
whose product is $w.$

Let $W_P$ denote the \textit{parabolic subgroup} of $W$ 
consisting of the element $w$ such that $w(\{1,\ldots,n\})\subset \{1,\ldots,n\}.$
Clearly $W_P$ is isomorphic to $S_n.$
Let $W^P$ denote the 
set of $w\in W$ such that 
$w(1)<\cdots<w(n).$
Let $u\in W.$ The coset $uW_P$ 
contains a unique element $w$ in $W^P.$
Actually $w$ is the unique 
element in the coset $uW_P$
of minimal length.
The longest element in $W$ is denoted by $w_0.$
If the coset $uW_P\in W/W_P$ is represented by $w\in W^P$,
then $w_0uW_P$ is represented by 
$w^\vee=(\overline{w(n)},\ldots,\overline{w(1)})\in W^P.$

\subsection{Combinatorial description of $W^P$}
By a \textit{symmetric} Young diagram, we mean
a sequence $D=(d_1,\ldots,d_n)$ 
with $d_1\geq d_2\geq \cdots\geq d_n\geq 0$ such that 
$$
d_i=\sharp\{j\;|\;d_j\geq i\}\quad (1\leq i\leq n).
$$
By $\mathcal{Y}^{sym}_n$ we denote 
the set of symmetric Young diagrams $D=(d_1,\ldots,d_n)$
contained in the square $n\times n,$
where by the last condition we mean $d_1\leq n.$

Let $w\in W^P.$ Then the sequence 
\beqn
D(w)=(n+1-w(1),n+2-w(2),\ldots,2n-w(n)),\label{D(w)}
\eeqn
is an element of $\mathcal{Y}^{sym}_n.$
Note that here we consider $w(i)$ simply  
as an element of $\{1,\ldots,2n\}$
without ``\textit{bar}''. 
For example if $n=5$ and $w=(1,3,4,\overline{5},\overline{2})
=(1,3,4,6,9),$ then
the corresponding Young diagram is $D(w)=(5,4,4,3,1).$
See Fig. 1.

\bigskip
$$\begin{array}{cc}
\unitlength 0.1in
\begin{picture}( 10.0000, 10.0000)(  2.0000,-12.0000)
%
\special{pn 13}%
\special{pa 200 200}%
\special{pa 400 200}%
\special{pa 400 400}%
\special{pa 200 400}%
\special{pa 200 200}%
\special{fp}%
%
\special{pn 13}%
\special{pa 400 200}%
\special{pa 600 200}%
\special{pa 600 400}%
\special{pa 400 400}%
\special{pa 400 200}%
\special{fp}%
%
\special{pn 13}%
\special{pa 400 400}%
\special{pa 600 400}%
\special{pa 600 600}%
\special{pa 400 600}%
\special{pa 400 400}%
\special{fp}%
%
\special{pn 13}%
\special{pa 600 200}%
\special{pa 800 200}%
\special{pa 800 400}%
\special{pa 600 400}%
\special{pa 600 200}%
\special{fp}%
%
\special{pn 13}%
\special{pa 600 400}%
\special{pa 800 400}%
\special{pa 800 600}%
\special{pa 600 600}%
\special{pa 600 400}%
\special{fp}%
%
\special{pn 13}%
\special{pa 600 600}%
\special{pa 800 600}%
\special{pa 800 800}%
\special{pa 600 800}%
\special{pa 600 600}%
\special{fp}%
%
\special{pn 13}%
\special{pa 800 200}%
\special{pa 1000 200}%
\special{pa 1000 400}%
\special{pa 800 400}%
\special{pa 800 200}%
\special{fp}%
%
\special{pn 13}%
\special{pa 800 400}%
\special{pa 1000 400}%
\special{pa 1000 600}%
\special{pa 800 600}%
\special{pa 800 400}%
\special{fp}%
%
\special{pn 13}%
\special{pa 800 600}%
\special{pa 1000 600}%
\special{pa 1000 800}%
\special{pa 800 800}%
\special{pa 800 600}%
\special{fp}%
%
\special{pn 13}%
\special{pa 1000 200}%
\special{pa 1200 200}%
\special{pa 1200 400}%
\special{pa 1000 400}%
\special{pa 1000 200}%
\special{fp}%
%
\special{pn 13}%
\special{pa 200 400}%
\special{pa 400 400}%
\special{pa 400 600}%
\special{pa 200 600}%
\special{pa 200 400}%
\special{fp}%
%
\special{pn 13}%
\special{pa 200 600}%
\special{pa 400 600}%
\special{pa 400 800}%
\special{pa 200 800}%
\special{pa 200 600}%
\special{fp}%
%
\special{pn 13}%
\special{pa 200 800}%
\special{pa 400 800}%
\special{pa 400 1000}%
\special{pa 200 1000}%
\special{pa 200 800}%
\special{fp}%
%
\special{pn 13}%
\special{pa 400 600}%
\special{pa 600 600}%
\special{pa 600 800}%
\special{pa 400 800}%
\special{pa 400 600}%
\special{fp}%
%
\special{pn 13}%
\special{pa 400 800}%
\special{pa 600 800}%
\special{pa 600 1000}%
\special{pa 400 1000}%
\special{pa 400 800}%
\special{fp}%
%
\special{pn 13}%
\special{pa 600 800}%
\special{pa 800 800}%
\special{pa 800 1000}%
\special{pa 600 1000}%
\special{pa 600 800}%
\special{fp}%
%
\special{pn 13}%
\special{pa 200 1000}%
\special{pa 400 1000}%
\special{pa 400 1200}%
\special{pa 200 1200}%
\special{pa 200 1000}%
\special{fp}%
\end{picture}%
\hspace{3cm}&
\unitlength 0.1in
\begin{picture}( 10.0000, 10.0000)(  2.0000,-12.0000)
%
\special{pn 13}%
\special{pa 200 200}%
\special{pa 400 200}%
\special{pa 400 400}%
\special{pa 200 400}%
\special{pa 200 200}%
\special{fp}%
%
\special{pn 13}%
\special{pa 400 200}%
\special{pa 600 200}%
\special{pa 600 400}%
\special{pa 400 400}%
\special{pa 400 200}%
\special{fp}%
%
\special{pn 13}%
\special{pa 400 400}%
\special{pa 600 400}%
\special{pa 600 600}%
\special{pa 400 600}%
\special{pa 400 400}%
\special{fp}%
%
\special{pn 13}%
\special{pa 600 200}%
\special{pa 800 200}%
\special{pa 800 400}%
\special{pa 600 400}%
\special{pa 600 200}%
\special{fp}%
%
\special{pn 13}%
\special{pa 600 400}%
\special{pa 800 400}%
\special{pa 800 600}%
\special{pa 600 600}%
\special{pa 600 400}%
\special{fp}%
%
\special{pn 13}%
\special{pa 600 600}%
\special{pa 800 600}%
\special{pa 800 800}%
\special{pa 600 800}%
\special{pa 600 600}%
\special{fp}%
%
\special{pn 13}%
\special{pa 800 200}%
\special{pa 1000 200}%
\special{pa 1000 400}%
\special{pa 800 400}%
\special{pa 800 200}%
\special{fp}%
%
\special{pn 13}%
\special{pa 800 400}%
\special{pa 1000 400}%
\special{pa 1000 600}%
\special{pa 800 600}%
\special{pa 800 400}%
\special{fp}%
%
\special{pn 13}%
\special{pa 800 600}%
\special{pa 1000 600}%
\special{pa 1000 800}%
\special{pa 800 800}%
\special{pa 800 600}%
\special{fp}%
%
\special{pn 13}%
\special{pa 1000 200}%
\special{pa 1200 200}%
\special{pa 1200 400}%
\special{pa 1000 400}%
\special{pa 1000 200}%
\special{fp}%
%
\special{pn 4}%
\special{pa 200 1200}%
\special{pa 200 400}%
\special{da 0.020}%
\special{pa 200 1200}%
\special{pa 400 1200}%
\special{da 0.020}%
\special{pa 400 1200}%
\special{pa 400 1000}%
\special{da 0.020}%
\special{pa 400 1000}%
\special{pa 800 1000}%
\special{da 0.020}%
\special{pa 800 1000}%
\special{pa 800 800}%
\special{da 0.020}%
%
\special{pn 4}%
\special{pa 400 600}%
\special{pa 400 1000}%
\special{da 0.020}%
\special{pa 600 800}%
\special{pa 200 800}%
\special{da 0.020}%
\special{pa 200 600}%
\special{pa 400 600}%
\special{da 0.020}%
\special{pa 200 1000}%
\special{pa 400 1000}%
\special{da 0.020}%
\special{pa 600 1000}%
\special{pa 600 800}%
\special{da 0.020}%
\end{picture}%
\\
\mbox{{Fig.1}}\hspace{3cm}&\mbox{{Fig.2}}
\end{array}
$$

For any symmetric Young diagram $D=(d_1,\ldots,d_n)$ in $\mathcal{Y}_n^{sym},$
its \textit{upper shifted diagram} $S(D)$
is obtained from $D$ by discarding the 
boxes strictly lower than the diagonal, i.e. 
$$
S(D)=\{(i,j)\in \Z^2\;|\; 1\leq i\leq j\leq d_i\},
$$
which we regard as an array of boxes in the plane
with matrix-style coordinates.
For example if $D=(5,4,4,3,1),$ its 
upper shifted diagram $S(D)$
is depicted as Fig.2 above.

Let $\la_i$ be the number of boxes in the $i$-th row
of $S(D).$ Then the sequence $\la=(\la_1,\ldots,\la_n)$ is
a \textit{strict partition}.
Namely there is $k$ such that
$\la_1>\cdots>\la_k>0$ and $\la_j=0$ for $j>k.$
Let $\mathcal{SP}_n$ denote the set of strict partitions
$\la=(\la_1,\ldots,\la_n)$ contained in the 
\textit{staircase} $\rho(n)=(n,n-1,\ldots,1),$
namely $\la_1\leq n.$
For $D$ in $\mathcal{Y}^{sym}_n,$ 
its upper shifted diagram $S(D)$ 
is thus considered to be a strict partition in $\mathcal{SP}_n.$
For example, the diagram of Fig. 2 is considered to be 
the strict partition $\la=(5,3,2).$

We let $\mathcal{M}_n$ denote the set $\{0,1\}^n.$
We use $\delta=(\delta_1,\ldots,\delta_n)$ 
to denote an element in $\mathcal{M}_n.$
For $w\in W^P$, we set $\delta_i=1$ if $i\in \{w(1),\ldots,w(n)\}$ and
$\delta_i=0$ if $i\notin \{w(1),\ldots,w(n)\}.$
Then we associate $\delta=(\delta_1,\ldots,\delta_n)
\in \mathcal{M}_n$ to $w\in W^P.$

\begin{prop} 
By the above correspondences,
we have bijections between the following sets: 

(i) The coset representatives $W^P$; 

(ii) The set $\mathcal{Y}_n^{sym}$ of symmetric Young diagrams
contained in the square $n\times n$; 

(iii) The set $\mathcal{M}_n$ of 
sequences $\delta=(\delta_1,\ldots,\delta_n)$ with $\delta_i\in\{0,1\}$; 

(iv) The set $\mathcal{SP}_n$ of strict partitions
$\la$ 
contained in $\rho(n).$
\end{prop}
\textit{Proof.} It is clear that each
correspondence is one to one.
Also it is easy to see that the cardinarity of each of the sets 
$W^P,\mathcal{Y}_n^{sym},\mathcal{M}_n,\mathcal{SP}_n$
is $2^n.$ Hence the claim follows.
$\square$

\bigskip
The following result is well known (see e.g. \cite{H}).

\begin{lem}\label{length}
For $w\in W^P,$ 
the length $\ell(w^\vee)$ is 
equal to $|\la|=\sum_{i=1}^k \la_i$,
where $\la\in \mathcal{SP}_n$ corresponds to $w\in W^P.$
\end{lem}

\section{Lagrangian Grassmannians and 
Schubert varieties}\label{Lag}
\setcounter{equation}{0}

This section is devoted to the set up of geometric objects. 
	\subsection{Lagrangian Grassmannians}

  Let $V$ be a vector space spanned by 
the basis $\pmb{e}_1,\dots,\pmb{e}_{n},\e_{\bar{n}},\ldots,\e_{\bar{1}}.$  
Introduce a symplectic form by 
$(\pmb{e}_i,\pmb{e}_j)=(\e_{\bar{i}},\e_{\bar{j}})=0$ and
$(\pmb{e}_i,\e_{\bar{j}})
=-(\e_{\bar{j}},\pmb{e}_i)=\delta_{ij}.$
Let $V_i$ denote the subspace spanned by the first $i$ 
vectors in $\pmb{e}_1,\dots,\pmb{e}_n,\e_{\bar{n}},\ldots,\e_{\bar{1}}.$
A subspace $W$ in $V$ is \textit{isotropic} if $(\pmb{u},\pmb{v})=0$ 
for all $\pmb{u},\pmb{v}\in W.$ 
Note that $V_i\;(1\leq i\leq n)$ are isotropic
of dimension $i$
and $V_{n+i}=(V_{n-i})^{\perp}\;(1\leq i\leq n).$
Denote by $LG_n$ the set of $n$-dimensional isotropic subspaces
of $V.$ Then $LG_n$ is a closed subvariety of 
the Grassmannian of $n$-dimensional subspaces of $V,$
and is called the \textit{Lagrangian Grassmannian}.

	The group $G=Sp(V)$ of linear automorphisms of $V$
preserving $(\,,\,)$ acts transitively on $LG_n.$  
So $LG_n$ is identified with the quotient of $G$ by the
stabilizer of any point. We identify $LG_n$ with the
homogeneous space $G/P$,
where $P$ denotes the stabilizer of the point $V_n$,
the span of $\e_1,\ldots,\e_n.$ 
The elements
of $G$ that are diagonal with respect to the
basis we took form a maximal
torus $T$ of $G.$ The elements of $G$ that
are upper triangular matrices form a Borel 
subgroup $B$ of $G.$

	\subsection{Schubert varieties}
The $T$-fixed points of $LG_n$ are parametrized
by $W^P:$ for  
$w$ in $W^P,$
the corresponding $T$-fixed point, denoted by $e(w),$
is the span of
$\e_{w(1)},\ldots,\e_{w(n)}.$
Let $X(w)^\circ$ denote the $B$-orbit
of $e(w).$ It is known 
thet $X(w)^\circ$ 
is an affine space of dimension $\ell(w),$
called a \textit{Schubert cell}.
The Zariski closures $X(w)=\overline{X(w)^\circ}$ are
called the \textit{Schubert varieties}.

We have the following description: 
$$
X(w)=\{L\in LG_n\;|\;\dim(L\cap V_{w(i)})\geq i\quad
\mbox{for}\;1\leq i\leq n\}.
$$
Let $\la=(\la_1>\cdots>\la_k>0)\in \mathcal{SP}_n$
be a strict partition 
corresponding to $w\in W^P.$
Then we also denote the variety $X(w)$ by $X(\la).$ 
We have
$$
X(\la)=\{
L\in LG_n\;|\;\dim(L\cap V_{n+1-\lambda_i})\geq i\quad
\mbox{for}\;1\leq i\leq k
\},
$$
whose codimension is given as $|\la|=\sum_{i=1}^k\la_i$ by Lemma \ref{length}.

	\subsection{Bruhat-Chevalley order}

For $w,v\in W^P$, we say $w\geq v$
if the torus fixed point ${e}(v)$ belongs to $X(w).$
This is a partial order called 
the \textit{Chevalley-Bruhat order.} 
The condition $w\geq v$ is given by the following 
two equivalent forms:
(1) $w(i)\geq v(i)\;(1\leq i\leq n),$
(2) $\la\subset \mu$, namely $\la_i\leq \mu_i\;(1\leq i\leq n),$
where $\la,\mu\in \mathcal{SP}_n$ correspond to $w,v$ respectively.
The Schubert variety 
$X(\la)$ admits a canonical partition
into Schubert cells $X(\mu)^\circ$ with $\mu\leq \la.$

\section{Equivariant cohomology}\label{Equ}
\setcounter{equation}{0}

We are interested in the $T$-equivariant 
integral cohomology ring $H_T^*(LG_n).$
For general facts on the equivariant cohomology,
we refer to Brion \cite{B},
Goresky, Kottwitz, and MacPherson \cite{GKM}
and references therein.

	\subsection{Equivariant Schubert classes}\label{ESC}
Let $\mathcal{S}$ denote the $T$-equivariant 
integral cohomology ring of a point $\{pt\}$
(namely the ordinary integral cohomology
ring of the classifying space of $T$).
The natural map $LG_n\rightarrow \{pt\}$ induce 
an $\mathcal{S}$-algebra structure on 
$H_T^*(LG_n).$
Given
$w\in W^P,$ denote by $\sigma(w)$ the $T$-equivariant fundamental class
of $X(w)$, 
called the 
\textit{equivariant Schubert class}.
We also denote $\sigma(w)$ by $\sigma(\la)$, where  
$\la\in \mathcal{SP}_n$ corresponds to $w\in W^P.$
It is known that  
$H_T^*(LG_n)$ is a 
free $\mathcal{S}$-module 
with the basis $\sigma(w)\;(w\in W^P):$
$$
H_T^*(LG_n)=\bigoplus_{w\in W^P}\mathcal{S}\cdot \sigma(w).
$$

For each $T$-fixed point $e(v),$ $v\in W^P$,
we have an embedding $i_v:\{e(v)\}\hookrightarrow LG_n.$
This yields a homomorphism 
$i_v^*:H_T^*(LG_n)\rightarrow H_T^*(\{e(v)\})\cong \mathcal{S}.$
The direct product of these 
is an injection of rings:
\begin{equation}
{\small{\prod_{v}}}i_v^*:
H_T^*(LG_n)\longrightarrow \prod_{v}H_T^*(\{e(v)\}).\label{inj}
\end{equation}
The injectivity is a consequence of 
``equivariant formality'' (cf. \cite{GKM}) of the $T$-variety $LG_n.$
For $w,v\in W^P$, denote by $\sigma(w)|_v$ the
image $i^*_v\sigma(w).$
The goal of this paper is 
to give an explicit
formula for $\sigma(w)|_v\in \mathcal{S}.$

\bigskip

\textbf{Remark.} A remarkable characterization of the image
of the morphism $\prod_{v}i_v^*$ has been obtained
by \cite{GKM}.
However we shall not use the result in the present paper.

\bigskip

Let $\mathfrak{t}=\mathrm{Lie}(T)$ be the Lie algebra of the torus $T.$
An element of $\mathfrak{t}$ takes the form
$h=\diag(h_1,\ldots,h_n,-h_n,\ldots,-h_1).$
Define linear functionals $\ep_i\in \mathfrak{t}^*$ by 
$\ep_i(h)=h_i$ for $1\leq i\leq n.$ 
Let $\widehat{T}$ be the free abelian group
generated by $\ep_1,\ldots,\ep_n.$ 
Each element $\sum m_i\ep_i$ in $\widehat{T}$ determines
a character of $T$ via
$T\ni \exp(h)\mapsto e^{\sum m_i\ep_i(h)}\in \C^\times.$
By this correspondence we can identify $\widehat{T}$ with
the character group of $T.$
There is a canonical map $\widehat{T}\rightarrow \mathcal{S}$
that extends to an isomorphism of the symmetric
algebra $\mathrm{Sym}(\widehat{T})$ 
onto $\mathcal{S}$ (see e.g. \cite{Br}, \S 1).
Thus we identify $\mathcal{S}$ with the polynomial ring  
$\Z[\ep_1,\ldots,\ep_n].$
We shall also use the variables $x_i=-\ep_i$ for $1\leq i\leq n.$
They are convenient for positivity reasons (cf. \S \ref{EPC}).

\subsection{$T$-stable affine neighborhood of $e(v)$}

Let $U(v)$ denote the set
of points of $LG_n$
whose matrix representatives $\xi=(\xi_{i,j})_{2n\times n}$ satisfy
$\xi_{v(i),j}=\delta_{i,j}.$
This is a $T$-stable affine space isomorphic to $\mathbb{A}^{n(n+1)/2}$
containing the point $e(v)$ as the origin.
The coordinate function  on $U(v)$ determined by
the matrix entry $\xi_{i,j}$ with $i\notin 
\{v(1),\ldots,v(n)\}$ 
is an eigenvector for $T$ of the 
weight $-(\ep_i-\ep_{v(j)}),$ here 
we understand $\ep_{\overline{i}}=-\ep_{i}.$
Let $\Xi$ be the square 
matrix $(\xi_{v^{\vee}(i),v(j)})_{1\leq i,j\leq n}.$
Not all the entries of $\Xi$ are independent,
since the column vectors of $\xi$ should span
an isotropic subspace. 
As a set of free parameters on the affine space $U(v)$, 
we can take the
set of entries of weakly upper triangular part of $\Xi$
with respect to the anti-diagonal.
Thus the coordinate ring
of $U(v)$ is $R(v)=\C[\xi_{v^{\vee}(i),v(j)}|v^{\vee}(i)\geq v(j)].$

It is convenient to
consider $v^\vee$ and $v$ as the (orderd) index sets
corresponding to the rows and columns of $\Xi$
respectively.
In this notation
we write $\Xi=(\xi_{r,c})_{r\in v^\vee,c\in v}.$
Note that the weight of the coordinate 
function $\xi_{r,c}$ is given by $-(\ep_r-\ep_c).$

\subsection{A product formula}

We shall prove a formula that expresses
$\sigma(w)|_w$ as a product 
of negative roots. 
This is 
a special case of the main result (Theorem \ref{res}).
In the subsequent of the paper, we need only 
the fact that $\sigma(w)|_w$ is 
a non-zero polynomial (see the Proof of Lemma  \ref{SpecialCoeff}).

\begin{lem}We have the following formula:
\beqn
\sigma(w)|_w=\prod_{(i,j)\in \la}(x_{w(i)}-x_{\overline{w(j)}}),
\label{diag}
\eeqn
where $\la\in \mathcal{SP}_n$ 
is the upper shifted diagram corresponding to $w\in W^P.$
\end{lem}

\textbf{Proof.} 
The variety $U(w)\cap X(w)$ is just a ``coordinate 
subspace'' in $U(w)$
defined by $\xi_{r,c}=0$ 
for $r\in w^\vee,\,c\in w,$ and $r>c.$
We denote by $\mathcal{I}(w)$ the set of
such pairs.  
For each $(r,c)\in\mathcal{I}(w)$, we associate $(i,j)$ 
by $w(i)=c$ and $\overline{w(j)}=r.$
Then $(i,j)$ is a box in the upper shifted diagram
$\la$ corresponding to $w.$
This establishes a bijection
from $\mathcal{I}(w)$ to the set of boxes of $\la.$
Recall that $\xi_{r,c}$ has 
the weight $-(\ep_r-\ep_c).$ 
Then by using Theorem 3 in \cite{LRS} 
we obtain the formula.
$\square$

\bigskip

For example,
let $w=(1,3,\bar{5},\bar{4},\bar{2}).$
The corresponding strict partition is $\la=(5,3).$
Then $\sigma(w)|_w=
2x_1(x_1+x_3)(x_1-x_5)(x_1-x_4)(x_1-x_2)
\times 2x_3(x_3-x_5)(x_3-x_4).
$

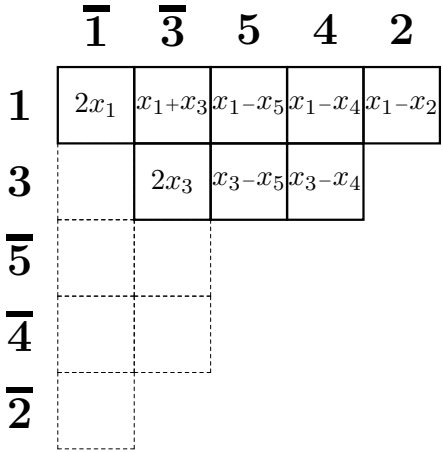
\begin{center}
\hspace{-6cm}
\begin{small}
\unitlength 0.1in
\begin{picture}( 49.9000, 23.1500)(-25.9000,-24.0000)
%
\special{pn 13}%
\special{pa 400 400}%
\special{pa 800 400}%
\special{pa 800 800}%
\special{pa 400 800}%
\special{pa 400 400}%
\special{fp}%
%
\special{pn 13}%
\special{pa 800 400}%
\special{pa 1200 400}%
\special{pa 1200 800}%
\special{pa 800 800}%
\special{pa 800 400}%
\special{fp}%
%
\special{pn 13}%
\special{pa 1200 400}%
\special{pa 1600 400}%
\special{pa 1600 800}%
\special{pa 1200 800}%
\special{pa 1200 400}%
\special{fp}%
%
\special{pn 13}%
\special{pa 1600 400}%
\special{pa 2000 400}%
\special{pa 2000 800}%
\special{pa 1600 800}%
\special{pa 1600 400}%
\special{fp}%
%
\special{pn 13}%
\special{pa 2000 400}%
\special{pa 2400 400}%
\special{pa 2400 800}%
\special{pa 2000 800}%
\special{pa 2000 400}%
\special{fp}%
%
\special{pn 13}%
\special{pa 800 800}%
\special{pa 1200 800}%
\special{pa 1200 1200}%
\special{pa 800 1200}%
\special{pa 800 800}%
\special{fp}%
%
\special{pn 13}%
\special{pa 1200 800}%
\special{pa 1600 800}%
\special{pa 1600 1200}%
\special{pa 1200 1200}%
\special{pa 1200 800}%
\special{fp}%
%
\special{pn 13}%
\special{pa 1600 800}%
\special{pa 2000 800}%
\special{pa 2000 1200}%
\special{pa 1600 1200}%
\special{pa 1600 800}%
\special{fp}%
%
\special{pn 4}%
\special{pa 800 800}%
\special{pa 400 800}%
\special{pa 400 1200}%
\special{pa 800 1200}%
\special{pa 800 800}%
\special{da 0.020}%
%
\special{pn 4}%
\special{pa 800 1200}%
\special{pa 400 1200}%
\special{pa 400 1600}%
\special{pa 800 1600}%
\special{pa 800 1200}%
\special{da 0.020}%
%
\special{pn 4}%
\special{pa 800 1600}%
\special{pa 400 1600}%
\special{pa 400 2000}%
\special{pa 800 2000}%
\special{pa 800 1600}%
\special{da 0.020}%
%
\special{pn 4}%
\special{pa 800 2000}%
\special{pa 400 2000}%
\special{pa 400 2400}%
\special{pa 800 2400}%
\special{pa 800 2000}%
\special{da 0.020}%
%
\special{pn 4}%
\special{pa 1200 1200}%
\special{pa 800 1200}%
\special{pa 800 1600}%
\special{pa 1200 1600}%
\special{pa 1200 1200}%
\special{da 0.020}%
%
\special{pn 4}%
\special{pa 1200 1600}%
\special{pa 800 1600}%
\special{pa 800 2000}%
\special{pa 1200 2000}%
\special{pa 1200 1600}%
\special{da 0.020}%
\put(6.0000,-6.0000){\makebox(0,0){${ 2{x}_1}$}}%
\put(10.0000,-6.0000){\makebox(0,0){${ {x}_1{\scriptscriptstyle +}{x}_3}$}}%
\put(14.0000,-6.0000){\makebox(0,0){${x}_1{\scriptscriptstyle -}{x}_5$}}%
\put(18.0000,-6.0000){\makebox(0,0){${x}_1{\scriptscriptstyle -}{x}_4$}}%
\put(22.0000,-6.0000){\makebox(0,0){${x}_1{\scriptscriptstyle -}{x}_2$}}%
\put(10.0000,-10.0000){\makebox(0,0){${ 2{x}_3}$}}%
\put(14.0000,-10.0000){\makebox(0,0){${x}_3{\scriptscriptstyle -}{x}_5$}}%
\put(18.0000,-10.0000){\makebox(0,0){${x}_3{\scriptscriptstyle -}{x}_4$}}%
\put(6.0000,-1.7000){\makebox(0,0){{\Large {${\displaystyle \above2pt \mbox{{\Large {\bf 1}}}}$}}}}%
\put(10.0000,-1.7000){\makebox(0,0){{\Large {${\displaystyle \above2pt \mbox{{\Large {\bf 3}}}}$}}}}%
\put(14.0000,-2.0000){\makebox(0,0){{\Large {\bf 5}}}}%
\put(18.0000,-2.0000){\makebox(0,0){{\Large {\bf 4}}}}%
\put(22.0000,-2.0000){\makebox(0,0){{\Large {\bf 2}}}}%
\put(2.0000,-6.0000){\makebox(0,0){{\Large {\bf 1}}}}%
\put(2.0000,-10.0000){\makebox(0,0){{\Large {\bf 3}}}}%
\put(2.0000,-13.7000){\makebox(0,0){{\Large {${\displaystyle \above2pt \mbox{{\Large {\bf 5}}}}$}}}}%
\put(2.0000,-17.7000){\makebox(0,0){{\Large {${\displaystyle \above2pt \mbox{{\Large {\bf 4}}}}$}}}}%
\put(2.0000,-21.7000){\makebox(0,0){{\Large {${\displaystyle \above2pt \mbox{{\Large {\bf 2}}}}$}}}}%
\end{picture}%
\end{small}
\end{center}

	\subsection{The divisor class}
	
Let $\mathrm{div}=(n,\overline{n-1},\ldots,\overline{2},\overline{1}).$ 
The corresponding $\sigma(\mathrm{div})$
is the unique Schubert class of codimension one. 
So we call it the \textit{divisor class}.
We know the following  
explicit form of this class
restricted to any $T$-fixed point $e(v).$ 
\begin{lem}\label{divclass}
The restriction of the divisor class $\sigma(\mathrm{div})$ to
a $T$-fixed point $e(v)$ is given by
\beqn
\sigma(\mathrm{div})|_{v}
=2\sum_{i=1}^n \delta_i x_i
\quad (v\in W^P),
\label{DivClass}
\eeqn
where $\delta=(\delta_1,\ldots,\delta_n)\in \mathcal{M}_n$ corresponds
to $v\in W^P.$
\end{lem}
\textit{Proof.}
Consider 
the closed subvariety $U(v)\cap X(\mathrm{div})$ of $U(v)\cong 
\mathbb{A}^{n(n+1)/2}.$
Let $\xi=[\pmb{\xi}_{v(1)},\ldots,\pmb{\xi}_{v(n)}]$ be 
a matrix represenatative of a point $L$ in $U(v).$ 
The condition for $L$ to be in $X(\mathrm{div})$ is 
equivalent to $\dim(V_n+L)\leq 2n-1.$
If we define the $n\times n$ matrix $X$ by
$$
\left[
\e_1,\ldots,\e_n,\pmb{\xi}_{v(1)},\ldots,\pmb{\xi}_{v(n)}
\right]
=\left[
\begin{array}{cc}
1_n & * \\
 0  & X
\end{array}
\right],
$$
then the last condition says
that $\det X$ 
should vanish. 
Let $k$ be such a number that $v(k)\leq n$ and $v(k+1)>n.$
Then by elementary
manipulations of 
a determinant, we have 
$\det X=
\pm \det Y$ where we denote by $Y=(\xi_{r,c})$ 
the $k\times k$ submatrix of $X$ with 
$c\in \{v(1),\ldots,v(k)\},$ and $
r\in \{\overline{v(k)},\ldots,\overline{v(1)}\}.$
Introducing a suitable monomial order on $R(v)$ such 
that the initial term of $\det Y$ is the product of 
anti-diagonal entries $\pm\prod_{i=1}^k \xi_{\overline{v(i)},v(i)}.$
Now applying Theorem 3 in \cite{LRS} we have
$\sigma(\mathrm{div})|_{v}=-\sum_{i=1}^k 2\ep_{v(i)}.$
Hence the claim follows.
$\square$

\bigskip

\textbf{Remark.} 
To the flag variety of the Kac-Moody groups,
Kostant and Kumar derived the corresponding 
formula (\cite{KK}, Prop. 4.24 (c), 
see also \cite{Ku} \S 11).

\subsection{Chevalley's multiplicities}\label{ChM}
	
	Let us recall the Chevalley multiplicites \cite{Ch}.
Let $w,w'\in W^P,$
such that $X(w')$ is a \textit{Schubert divisor} of $X(w),$
i.e., $X(w')$ is a codimension one subvariety 
in $X(w).$
Then there is a positive root $\beta$ such that 
$w'=ws_{\beta}$ and $\ell(w')=\ell(w)-1,$ 
where $s_\beta$ is the reflection corresponding to $\beta.$
Let $(\,,\,)$ be the inner product on 
$\widehat{T}\otimes_\Z{\mathbb{R}}=\oplus_{i=1}^n\mathbb{R}\ep_i$
such that $(\ep_i,\ep_j)=\delta_{ij},$
and $\beta^\vee$ be $2\beta/(\beta,\beta).$
Then the Chevalley multiplicity $c(w,w')$ is 
defined 
\beqn
c(w,w')=( \varpi_n,\beta^{\vee}),
\eeqn
where $\varpi_n=\sum_{i=1}^n\ep_i,$ the $n$-th
fundamental weight.
We can describe $c(w,w)$ in a combinatorial way.

\begin{lem}(\cite{HB}) \label{CombinCh}
Let $X(w')$ be a Schubert divisor in $X(w).$
Let $D(w), D(w')$ be the
corresponding symmetric diagrams.
Exactly one of the following holds.

(1) $D(w')$ is obtained from $D(w)$ adding
two boxes at the positions $(i,j)$ and $(j,i)$ ($i\ne j$).
Then the corresponding positive root $\beta$ is ${\ep_i+\ep_j},$
and we have $c(w,w')=2.$ 

(2) $D(w')$ is obtained from $D(w)$ adding a box
at the diagonal position $(i,i)$, then the corresponding
positive root $\beta$ is $2\ep_i,$
and we have $c(w,w')=1.$
\end{lem}

In the following figure, the numbers 
indicate the Chevalley 
multiplicities, where $n=3$:
	\begin{center}
\unitlength 0.1in
\begin{picture}( 16.7000, 36.5500)( -2.6500,-39.8000)
%
\special{pn 13}%
\special{pa 120 1802}%
\special{pa 232 1802}%
\special{pa 232 1912}%
\special{pa 120 1912}%
\special{pa 120 1802}%
\special{fp}%
%
\special{pn 13}%
\special{pa 232 1912}%
\special{pa 344 1912}%
\special{pa 344 1802}%
\special{pa 232 1802}%
\special{pa 232 1912}%
\special{fp}%
%
\special{pn 13}%
\special{pa 344 1802}%
\special{pa 456 1802}%
\special{pa 456 1912}%
\special{pa 344 1912}%
\special{pa 344 1802}%
\special{fp}%
%
\special{pn 4}%
\special{pa 120 1912}%
\special{pa 232 1912}%
\special{pa 232 2024}%
\special{pa 120 2024}%
\special{pa 120 1912}%
\special{da 0.020}%
%
\special{pn 4}%
\special{pa 120 2024}%
\special{pa 232 2024}%
\special{pa 232 2136}%
\special{pa 120 2136}%
\special{pa 120 2024}%
\special{da 0.020}%
%
\special{pn 4}%
\special{pa 680 1522}%
\special{pa 790 1522}%
\special{pa 790 1410}%
\special{pa 680 1410}%
\special{pa 680 1522}%
\special{da 0.020}%
%
\special{pn 13}%
\special{pa 680 1298}%
\special{pa 790 1298}%
\special{pa 790 1410}%
\special{pa 680 1410}%
\special{pa 680 1298}%
\special{fp}%
%
\special{pn 13}%
\special{pa 790 1410}%
\special{pa 902 1410}%
\special{pa 902 1298}%
\special{pa 790 1298}%
\special{pa 790 1410}%
\special{fp}%
%
\special{pn 13}%
\special{pa 1182 1802}%
\special{pa 1294 1802}%
\special{pa 1294 1912}%
\special{pa 1182 1912}%
\special{pa 1182 1802}%
\special{fp}%
%
\special{pn 13}%
\special{pa 1294 1912}%
\special{pa 1406 1912}%
\special{pa 1406 1802}%
\special{pa 1294 1802}%
\special{pa 1294 1912}%
\special{fp}%
%
\special{pn 13}%
\special{pa 1406 1912}%
\special{pa 1294 1912}%
\special{pa 1294 2024}%
\special{pa 1406 2024}%
\special{pa 1406 1912}%
\special{fp}%
%
\special{pn 4}%
\special{pa 1294 2024}%
\special{pa 1182 2024}%
\special{pa 1182 1912}%
\special{pa 1294 1912}%
\special{pa 1294 2024}%
\special{da 0.020}%
%
\special{pn 13}%
\special{pa 624 2304}%
\special{pa 736 2304}%
\special{pa 736 2416}%
\special{pa 624 2416}%
\special{pa 624 2304}%
\special{fp}%
%
\special{pn 13}%
\special{pa 736 2416}%
\special{pa 846 2416}%
\special{pa 846 2304}%
\special{pa 736 2304}%
\special{pa 736 2416}%
\special{fp}%
%
\special{pn 13}%
\special{pa 846 2304}%
\special{pa 958 2304}%
\special{pa 958 2416}%
\special{pa 846 2416}%
\special{pa 846 2304}%
\special{fp}%
%
\special{pn 13}%
\special{pa 846 2528}%
\special{pa 736 2528}%
\special{pa 736 2416}%
\special{pa 846 2416}%
\special{pa 846 2528}%
\special{fp}%
%
\special{pn 4}%
\special{pa 624 2416}%
\special{pa 736 2416}%
\special{pa 736 2528}%
\special{pa 624 2528}%
\special{pa 624 2416}%
\special{da 0.020}%
%
\special{pn 4}%
\special{pa 624 2528}%
\special{pa 736 2528}%
\special{pa 736 2640}%
\special{pa 624 2640}%
\special{pa 624 2528}%
\special{da 0.020}%
%
\special{pn 13}%
\special{pa 624 2974}%
\special{pa 736 2974}%
\special{pa 736 3086}%
\special{pa 624 3086}%
\special{pa 624 2974}%
\special{fp}%
%
\special{pn 13}%
\special{pa 736 3086}%
\special{pa 846 3086}%
\special{pa 846 2974}%
\special{pa 736 2974}%
\special{pa 736 3086}%
\special{fp}%
%
\special{pn 13}%
\special{pa 846 2974}%
\special{pa 958 2974}%
\special{pa 958 3086}%
\special{pa 846 3086}%
\special{pa 846 2974}%
\special{fp}%
%
\special{pn 13}%
\special{pa 846 3198}%
\special{pa 736 3198}%
\special{pa 736 3086}%
\special{pa 846 3086}%
\special{pa 846 3198}%
\special{fp}%
%
\special{pn 4}%
\special{pa 624 3086}%
\special{pa 736 3086}%
\special{pa 736 3198}%
\special{pa 624 3198}%
\special{pa 624 3086}%
\special{da 0.020}%
%
\special{pn 4}%
\special{pa 624 3198}%
\special{pa 736 3198}%
\special{pa 736 3310}%
\special{pa 624 3310}%
\special{pa 624 3198}%
\special{da 0.020}%
%
\special{pn 4}%
\special{pa 736 3310}%
\special{pa 846 3310}%
\special{pa 846 3198}%
\special{pa 736 3198}%
\special{pa 736 3310}%
\special{da 0.020}%
%
\special{pn 13}%
\special{pa 846 3086}%
\special{pa 958 3086}%
\special{pa 958 3198}%
\special{pa 846 3198}%
\special{pa 846 3086}%
\special{fp}%
%
\special{pn 13}%
\special{pa 624 3646}%
\special{pa 736 3646}%
\special{pa 736 3758}%
\special{pa 624 3758}%
\special{pa 624 3646}%
\special{fp}%
%
\special{pn 13}%
\special{pa 736 3758}%
\special{pa 846 3758}%
\special{pa 846 3646}%
\special{pa 736 3646}%
\special{pa 736 3758}%
\special{fp}%
%
\special{pn 13}%
\special{pa 846 3646}%
\special{pa 958 3646}%
\special{pa 958 3758}%
\special{pa 846 3758}%
\special{pa 846 3646}%
\special{fp}%
%
\special{pn 13}%
\special{pa 846 3868}%
\special{pa 736 3868}%
\special{pa 736 3758}%
\special{pa 846 3758}%
\special{pa 846 3868}%
\special{fp}%
%
\special{pn 4}%
\special{pa 624 3758}%
\special{pa 736 3758}%
\special{pa 736 3868}%
\special{pa 624 3868}%
\special{pa 624 3758}%
\special{da 0.020}%
%
\special{pn 4}%
\special{pa 624 3868}%
\special{pa 736 3868}%
\special{pa 736 3980}%
\special{pa 624 3980}%
\special{pa 624 3868}%
\special{da 0.020}%
%
\special{pn 4}%
\special{pa 736 3980}%
\special{pa 846 3980}%
\special{pa 846 3868}%
\special{pa 736 3868}%
\special{pa 736 3980}%
\special{da 0.020}%
%
\special{pn 13}%
\special{pa 846 3758}%
\special{pa 958 3758}%
\special{pa 958 3868}%
\special{pa 846 3868}%
\special{pa 846 3758}%
\special{fp}%
%
\special{pn 13}%
\special{pa 846 3980}%
\special{pa 958 3980}%
\special{pa 958 3868}%
\special{pa 846 3868}%
\special{pa 846 3980}%
\special{fp}%
\put(7.9000,-4.1000){\makebox(0,0){{\Large$\phi$}}}%
%
\special{pn 13}%
\special{pa 846 852}%
\special{pa 736 852}%
\special{pa 736 962}%
\special{pa 846 962}%
\special{pa 846 852}%
\special{fp}%
\put(9.0200,-6.8300){\makebox(0,0){{\Large {\bf 1}}}}%
\put(9.0200,-11.3000){\makebox(0,0){{\Large {\bf 2}}}}%
\put(11.2600,-15.7700){\makebox(0,0){{\Large {\bf 1}}}}%
\put(4.5500,-15.7700){\makebox(0,0){{\Large {\bf 2}}}}%
\put(11.2600,-22.4800){\makebox(0,0){{\Large {\bf 2}}}}%
\put(4.5500,-22.4800){\makebox(0,0){{\Large {\bf 1}}}}%
\put(9.0200,-34.7700){\makebox(0,0){{\Large {\bf 1}}}}%
\put(9.0200,-28.0700){\makebox(0,0){{\Large {\bf 2}}}}%
%
\special{pn 8}%
\special{pa 790 796}%
\special{pa 790 572}%
\special{fp}%
\special{sh 1}%
\special{pa 790 572}%
\special{pa 770 638}%
\special{pa 790 624}%
\special{pa 810 638}%
\special{pa 790 572}%
\special{fp}%
%
\special{pn 8}%
\special{pa 790 1242}%
\special{pa 790 1018}%
\special{fp}%
\special{sh 1}%
\special{pa 790 1018}%
\special{pa 770 1086}%
\special{pa 790 1072}%
\special{pa 810 1086}%
\special{pa 790 1018}%
\special{fp}%
%
\special{pn 8}%
\special{pa 1126 1746}%
\special{pa 958 1578}%
\special{fp}%
\special{sh 1}%
\special{pa 958 1578}%
\special{pa 992 1638}%
\special{pa 996 1616}%
\special{pa 1020 1610}%
\special{pa 958 1578}%
\special{fp}%
%
\special{pn 8}%
\special{pa 958 2248}%
\special{pa 1126 2080}%
\special{fp}%
\special{sh 1}%
\special{pa 1126 2080}%
\special{pa 1066 2114}%
\special{pa 1088 2118}%
\special{pa 1094 2142}%
\special{pa 1126 2080}%
\special{fp}%
%
\special{pn 8}%
\special{pa 624 2248}%
\special{pa 456 2080}%
\special{fp}%
\special{sh 1}%
\special{pa 456 2080}%
\special{pa 488 2142}%
\special{pa 494 2118}%
\special{pa 516 2114}%
\special{pa 456 2080}%
\special{fp}%
%
\special{pn 8}%
\special{pa 790 3596}%
\special{pa 790 3372}%
\special{fp}%
\special{sh 1}%
\special{pa 790 3372}%
\special{pa 770 3438}%
\special{pa 790 3424}%
\special{pa 810 3438}%
\special{pa 790 3372}%
\special{fp}%
\special{pa 790 2924}%
\special{pa 790 2702}%
\special{fp}%
\special{sh 1}%
\special{pa 790 2702}%
\special{pa 770 2768}%
\special{pa 790 2754}%
\special{pa 810 2768}%
\special{pa 790 2702}%
\special{fp}%
%
\special{pn 8}%
\special{pa 456 1746}%
\special{pa 624 1578}%
\special{fp}%
\special{sh 1}%
\special{pa 624 1578}%
\special{pa 562 1610}%
\special{pa 586 1616}%
\special{pa 590 1638}%
\special{pa 624 1578}%
\special{fp}%
\end{picture}%
	\end{center}

One can easily verify the following rule:

\begin{lem}\label{lemlength} We assume $w'\to w\;(w,w'\in W^P).$
Let $\la,\la'\in \mathcal{SP}_n$ be correspond to $w,w'$ respectively.
Let $k,k'$ be the numbers of non-zero parts of 
$\la,\la'$ respectively. 
Then $c(w,w')=1$ if $k'=k+1$, and $c(w,w')=2$ if $k'=k.$
\end{lem}
	\subsection{Equivariant Pieri-Chevalley formula}\label{EPC}

Since $\{\sigma(w)\}_{w\in W^P}$ forms a basis
of $H_T^*(LG_n)$ over the ring $\mathcal{S}$, we can 
define the structure constants $c_{w,v}^u\in \mathcal{S}$ 
for all $w,v,u\in W^P$
 by the formula
\beqn
\sigma(w)\cdot\sigma(v)=\sum_{u}c^u_{w,v}\,\sigma(u).\label{EqCh}
\eeqn
The structure constant $c_{w,v}^u$ has degree $\ell(w)+\ell(v)-\ell(u)$ and 
vanishes unless $u\geq w,v$ and $\ell(u)\leq \ell(w)+\ell(v).$
It should be remarked that $c_{w,v}^u$ has
a remarkable positivity property 
conjectured by D. Peterson and proved by Graham \cite{G}.
Namely each $c_{w,v}^u$ can be written as a linear combination of monomials in 
the negative roots with 
nonnegative integer coefficients.

\begin{lem}\label{EPCformula}
(The equivariant Pieri-Chevalley formula) Let $w\in W^P.$
Then the following formula holds:
\beqn
\sigma(\mathrm{div})\cdot\sigma(w)
=c_{\mathrm{div},w}^w\;\sigma(w)
+\sum_{w':w'\to w}c(w,w')\,\sigma(w'),\label{Ch}
\eeqn
where $w'\to w$ means that $X(w')$ is a Schubert
divisor of $X(w).$ 
\end{lem}
\textit{Proof.} By the same argument
of \cite{KT}, Prop. 2, the claim follows from
the Pieri-Chevalley-type formulas 
for the ordinary integral cohomology, which were proved by Fulton and 
Woodward \cite{FW}, Lemma 8.1.
$\square$

\bigskip

For the flag variety of an arbitrary Kac-Moody group,
the corresponding formula of Lemma \ref{EPCformula}
has been appeared in 
the context of the ni-Hecke algebra
by Kostant and Kumar \cite{KK}.
Later
Arabia \cite{A} established the fact that 
the equivariant cohomology 
is isomorphic to the dual of the nil-Hecke algebra.
The parabolic analogue is also
studied in \cite{Ku}, \S 11. See also
Robinson \cite{R}, Andersen, Jantzen, and Sorgel \cite{AJS}, Appendix D.

	\subsection{Recurrence relation}

In this section, we prove a Key lemma (Lemma \ref{lemrec}) to the proof of 
our main result. First
we need a simple lemma on structure constants.

\begin{lem}
The structure constant $c_{\mathrm{div},w}^w$
is given by
\beqn
c_{\mathrm{div},w}^w=\sigma(\mathrm{div})|_{w}.\label{SpecialCoeff}
\eeqn
\end{lem}
\textit{Proof.}
If we restrict (\ref{Ch}) to $e(w)$, we have
$$
\sigma(\mathrm{div})|_{w}\cdot\sigma(w)|_{w}
=c_{\mathrm{div},w}^w\;\sigma(w)|_{w}
+\sum_{w':w'\to w}c(w,w')\,\sigma(w')|_{w}.
$$
For $w'$ such that $w'\to w,$ we have $\sigma(w')|_{w}=0$
since $w\not\leq w'.$
Hence the sum
in the right hand side 
vanishes.
The claim follows since $\sigma(w)|_{w}$ is non-zero
as we see from Lemma \ref{diag}. 
$\square$

\bigskip

Now the equivariant Pieri-Chevalley formula (\ref{Ch}) gives 
directly the following recurrence
relation on the family of restricted classes $\sigma(w)|_v\;(w\in W^P)$
for any fixed $v\in W^P.$

\begin{lem}\label{lemrec} Let $e(v)$ be any $T$-fixed point.
The polynomials $\sigma(w)|_v\,(w\in W^P)$ 
satisfy the
 following recurrence relation:
\beqn
d(w,v)\cdot\sigma(w)|_v
=\sum_{w':w'\rightarrow w}c(w,w')\, \sigma(w')|_v,\label{Rec}
\eeqn
where $d(w,v)=\sigma(\mathrm{div})|_{v}-\sigma(\mathrm{div})|_{w}.$
\end{lem}

Since $d(w,v)$ is non-zero if $w\geq v$ and $w\ne v,$
the recurrence relation (\ref{Rec})
and the initial condition $\sigma(\phi)|_v=1$
determine the polynomials
$\sigma(w)|_v\;(w\in W^P)$ uniquely.
An analogous recurrence relation was
used by Rosenthal and Zelevinsky \cite{RZ} to 
prove a determinantal formula of 
the multiplicity of a $T$-fixed point in
a Schubert variety in the Grassmannian. 

\bigskip

\textbf{Remark.} 
The ordinary-cohomology version of Lemma \ref{lemrec}
has been obtained by 
Lakshmibai and Weyman
\cite{LW}, and Hiller \cite{H}.

\section{The factorial Schur $Q$-functions}\label{Fac}
\setcounter{equation}{0}

Let $x=(x_1,\ldots,x_n)$ be a finite sequence of variables and 
let $a=(a_i)_{i\geq 1}$ be any sequence such that $a_1=0.$
Let $$
(x|a)^k=\prod_{i=1}^k(x-a_i)
$$
for each $k\geq 1$ and $(x|a)^0=1.$
The \textit{factorial} \textit{Schur} $Q$-\textit{function} 
for a strict partition $\la=(\la_1>\cdots>\la_k>0)$ of \textit{length} $k\leq n$ 
is defined
as follows \cite{Iv2}.
\begin{Def}\label{nimmo}
Let $A(x)$ denote the skew-symmetric
$n\times n$ matrix
$
({(x_i-x_j)}/{(x_i+x_j)})_{1\leq i,j\leq n}
$
and let $B_\la(x|a)$ denote the $n\times k$ matrix
$
((x_i|a)^{\la_{k-j+1}}).
$
Let
$$
A_\la(x|a)=
\left[
\begin{array}{cc}
A(x) & B_\la(x|a) \\
-{}^t\!B_\la(x|a) &0
\end{array}
\right]
$$
which is a skew-symmetric $(n+k)\times(n+k)$ matrix.
Put 
$$
\Pf_\la(x|a)=\begin{cases}
\Pf\left(A_\la(x_1,\ldots,x_n|a)\right)\quad 
\mbox{if} \;\; n+k\;\mbox{is even;}\\
\Pf\left(A_\la(x_1,\ldots,x_n,0\,|a)\right)\quad 
\mbox{if} \;\; n+k \;\mbox{is odd.}
\end{cases}
$$
Then put
\beqn
P_\la(x|a)=
\frac{\Pf_\la(x|a)}
{D_n(x)},\quad
Q_\la(x|a)=2^k P_\la(x|a),\label{Q}
\eeqn
where 
$
D_n(x)=\prod_{1\leq i<j\leq n}{(x_i-x_j)}/{(x_i+x_j)}.
$
\end{Def}

\textbf{Remark.}
The above definition is 
a factorial analogue of 
Nimmo's formula \cite{N} (see also
\cite{M1}, Ch. III, 8, Example 13) for the Schur $Q$-functions.
The reader can find other expressions for 
$Q_\la(x|a)$ in \cite{Iv2}.

\bigskip

The functions $Q_{\la}(x|a)$
were introduced by Ivanov\footnote{According to Ivanov \cite{Iv1,Iv2},
A.Okounkov defined them for
the special parameter $a$ with $a_i=i-1.$} \cite{Iv1, Iv2}.
He established some fundamental properties
of the functions (combinatorial presentations,
Schur-type Pfaffian formulas,
vanishing and characterization properties etc.).
In particular, a Pieri-type formula
is available, which is crucial to our consideration.
Note that 
$P_{(1)}(x|a)$ does not depend on the parameter $a=(a_i)$ and 
actually we have $P_{(1)}(x|a)=\sum_{i=1}^n{x_i}.$
So we simply denote $P_{(1)}(x|a)$ by $P_{(1)}(x).$
Let $\lambda$ and $\lambda'$ be 
strict partitions of length $\leq n.$
We will write $\lambda'\rightarrow \lambda$ if 
$\lambda\subset \lambda'$ and $|\lambda'|=|\lambda|+1.$

\begin{prop} \cite{Iv2} (A Pieri-type formula)
For any strict partition $\la=(\la_1>\cdots>\la_k>0)$ 
of length $k\leq n$, we have
\beqn
\left(P_{(1)}(x)-\sum_{j=1}^{k}a_{\la_j+1}\right)\cdot P_{\la}(x|a)
=\sum_{\la':\la'\to \la}P_{\la'}(x|a),\label{Pieri}
\eeqn
where $\la'$ runs for all strict partitions of length less than or
equal to $n$
such that $\lambda'\to\lambda.$
\end{prop}

{\bf Proof.} The above formula
corresponds to Thorem 6.2 in \cite{Iv2}.
The only difference is that we use $n$-variables
$x=(x_1,\ldots,x_n)$ here. Then
we can consistently set $P_\mu(x|a)$ to be zero
for any strict partition $\mu$ of length strictly greater than $n$
(see \cite{Iv2}, Definition 2.10). 
$\square.$

\bigskip

Factorial analogues of 
the  Schur $S$-functions were introduced by 
Biedenharn and Louck \cite{BL}
and further studied by Chen and Louck \cite{CL},
Goulden and Greene \cite{GG},
Goulden and Hamel \cite{GH}, Macdonald \cite{M2},
and Molev and Sagan \cite{MS}
(see also Macdonald \cite{M1}, Ch. I, 3, Example 20-21). 
In these works it was shown that several important
facts about the Schur $S$-functions (combinatorial presentations, 
Jacobi-Trudi identities, Pieri-type formulas, Littlewood-Richetchardson rules etc.) can be 
transferred to the factorial Schur $S$-functions.
The factorial Schur $S$-functions
also play a central role in the study
of the center of the universal enveloping 
algebra of $\mathfrak{gl}_n$ (see Okounkov and Olshanski 
\cite{OO}, Okounkov 
\cite{O} and references therein). 

In a geometric context, the factorial Schur
functions appeared in \cite{KT}, \cite{LRS}.
They present the restriction to torus fixed points
of the Schubert
classes in the equivariant cohomology 
of the Grassmannian.
Recently Mihalcea \cite{Mi} obtained 
a presentation by generators and relations
for the equivariant \textit{quantum} cohomology ring
of the Grassmannian. In this work
the factorial Schur $S$-functions\footnote{It has come to 
my knowledge  
via \cite{Mi} (\S5, Remark 2) that
the factorial Schur $S$-functions
coincide with 
the double Schubert polynomials
by Lascoux and Sch\"utzenberger \cite{LS}
when indexed by a Grassmann permutation.
However, the details of this connection seems to be missing
from the literature.} 
appeared as the polynomial representatives
of the equivariant quantum Schubert classes. 
A similar presentation
for the quantum cohomology ring
of the Lagrangian Grassmannian 
was given by Kresch and Tamvakis \cite{KrT}.
It will be an interesting problem
to extend their result 
to the quantum equivariant cohomology ring.

\bigskip

\section{Restriction and Giambelli-type formulas}\label{Main}
\setcounter{equation}{0}
\subsection{Restriction formula}\label{Res}
Let us take the following particular
parameters: 
\beqn
a_1=0,\quad
a_i=x_{n-i+2}\quad(2\leq i\leq n+1),\quad
a_i=0\quad (i>n+1).\label{a}
\eeqn
We denote 
by $x_{\langle n\rangle}$ the specialization 
of $a=(a_i)_{i\geq 1}$ given by (\ref{a}).
Let $\mu\in W^P$
and $\delta=(\delta_1,\ldots,\delta_n)\in \mathcal{M}_n$
be correspond to $\mu.$
Then we set 
$x_\mu=(\delta_1x_1,\ldots,\delta_nx_n).
$
\begin{Def}\label{specialization}
A specialization $Q_\lambda(x_\mu|x_{\langle n\rangle})$ of $Q_\lambda(x|a)$
is given as follows:
First we substitute $x_\mu$ for $x=(x_1,\ldots,x_n)$ 
to obtain $Q_\lambda(x_\mu|a)$, then 
we specialize $a_i$'s as in (\ref{a}) to get $Q_\lambda(x_\mu|x_{\langle n\rangle}).$
\end{Def}

\begin{thm}\label{res} For strict partitions $\la,\mu\in \mathcal{SP}_n,$ 
we have
\beqn
\sigma(\la)|_\mu=
Q_\la(x_\mu|x_{\langle n\rangle}).\label{MAIN}
\eeqn
\end{thm}

\textit{Proof.} 
It suffices to show that the right hand side 
of (\ref{MAIN}) satisfy 
the recurrence relation (\ref{Rec}) 
and the initial condition $\sigma(\phi)|_\mu=1.$
As for the initial condition, we have $Q_\phi(x|a)=1$
by definition.
Hence the proof is completed by a comparison
of (\ref{Rec}) and (\ref{Pieri}).
We first specializes $P_\lambda(x|a)$ to $P_\lambda(x_\mu|a)$
and then to $P_\lambda(x_\mu|x_{\langle n\rangle}).$
By applying this specialization to (\ref{Pieri}), we have
\beqn
\left(P_{(1)}(x_\mu)-\sum_{j=1}^{k}x_{n-\la_j+1}\right)\cdot P_{\la}(x_\mu|x_{\langle n\rangle})
=\sum_{\la'}P_{\la'}(x_\mu|x_{\langle n\rangle}),\label{Pieri_sp}
\eeqn
where the sum 
is taken over those $\lambda'\in\mathcal{SP}_n$ such that $\lambda'\to\lambda$
because $P_{\la'}(x_\mu|x_{\langle n\rangle})$ vanishes
unless $\la'\in \mathcal{SP}_n$ (Proposition \ref{vanish}).
Now we multiply the both hand sides of 
(\ref{Pieri_sp}) by $2^{k+1},$ where $k$ is the number of non-zero
parts of $\la.$
By Lemma \ref{divclass}, we have
$\sigma(\mathrm{div})|_w=2\sum_{j=1}^{k}x_{n-\la_j+1}$
and $\sigma(\mathrm{div})|_v=2P_{(1)}(x_\mu).$
Therefore we have
$$d(w,v)=2P_{(1)}(x_\mu)
-2\sum_{i=1}^k x_{n-\la_i+1}.$$
Now let $\la'\in \mathcal{SP}_n$ be such that $\la'\to \la$ and $k'$ be the number of non-zero
parts of $\la'.$
From Lemmas \ref{CombinCh} and \ref{lemlength},
we can see that $2^{k+1}P_{\la'}(x|a)=c(w,w')Q_{\la'}(x|a).$
Note also $2^k P_\la(x|a)=Q_\la(x|a).$
Thus we proved that $Q_\la(x_\mu|x_{\langle n\rangle})\;(\la\in \mathcal{SP}_n)$
satisfy (\ref{Rec}).
$\square$ 

\subsection{Giambelli-type formula}
Now we can prove
an equivariant analogue of Pragacz' 
Giambelli-type formula.
Let $\lambda\in \mathcal{SP}_n.$
We write $\lambda=(\lambda_1,\ldots,\lambda_{2r})$
with $\lambda_1>\cdots>\lambda_{2r}\geq 0.$
\begin{thm}\label{Giam}
The equivariant Schubert class $\sigma(\lambda)$ 
is expressed as a Pfaffian 
of the following form
\beqn
\sigma(\lambda)=
\mathrm{Pf}\left(
\sigma(\lambda_i,\lambda_j)
\right)_{1\leq i,j\leq 2r}.\label{Giambelli}
\eeqn
\end{thm}

\textbf{Proof.} 
Because of the injection (\ref{inj}),
it is enough to show that, for 
arbitrary $\mu$ in $\mathcal{SP}_n$,
the restrictions to $e(\mu)$
of the both hand sides
of (\ref{Giambelli}) coincide. 
We have
$$
\sigma(\lambda)|_{\mu}=Q_{\lambda}(x_{\mu}|x_{\langle n\rangle})
=
\mathrm{Pf}\left(
Q_{\lambda_i,\lambda_j}(x_{\mu}|x_{\langle n\rangle})
\right)_{1\leq i,j\leq 2r}
=
\mathrm{Pf}\left(
\sigma(\lambda_i,\lambda_j)|_{\mu}
\right)_{1\leq i,j\leq 2r}.
$$
In the second equality,
we use the Pfaffian formula 
for factorial $Q$-functions (\ref{IvSchur}).
Since the restriction $i_\mu^*$ is a ring 
homomorphism and we are done.
$\square.$

\bigskip

The above formula
has a striking character 
in contrast to the ordinary Grassmannian case \cite{LRS},
where the equivariant Giambelli formula
is given in a Jacobi-Trudi type determinant,
with matrix entries of {\it linear combinations} of (equivariant ) special Schubert classes. 
In our formula, each matrix entry of the Pfaffian 
is itself an equivariant Schubert class.
In spite of this simplicity, if we wish to
express the equivariant Schubert class as 
a polynomial of the {\it special Schubert classes} $\sigma(k)\;(1\leq k\leq n)$,
we need some work to be done. 
We will treat the problem in the next subsection. 

\section{On the two-row type classes}\label{Two}
\setcounter{equation}{0}
The formula (\ref{Giambelli}) looks 
the same as the classical 
one shown by Pragacz (\cite{P}, Prop. 6.6, see also J\'ozefiak \cite{Jo}), where
the $Q$-functions $Q_\lambda(x)\;(\lambda\in \mathcal{SP}_n)$ 
represent
the Schubert classes in the ordinary cohomology ring of $LG_n.$
Recall that we have the following formula for $r> s\geq 0$:
\beqn
Q_{r,s}(x)=Q_r(x)Q_s(x)+2\sum_{i=1}^s(-1)^iQ_{r+i}(x)Q_{s-i}(x).\label{quad}
\eeqn
Therefore the Pragacz' formula
gives an expression for each Schubert class
as a polynomial in the special Schubert classes.

Now in our setting of equivariant cohomology ring, 
equation (\ref{Giambelli}) actually provides an expression for each $\sigma(\lambda)$ as 
a polynomial in $\sigma(\lambda_i,\lambda_j).$
If $\lambda_j=0$ then $\sigma(\lambda_i,\lambda_j)=\sigma(\lambda_i)$ is a
{\it special\/} class.
For the {\it two-row} type classes, i.e. $\sigma(\lambda_i,\lambda_j)$ with $\lambda_j>0$,
we want to express them as a polynomial in
the special classes $\sigma(k)\;(1\leq k\leq n).$
In fact, we have the following 
expression for two-row type classes $\sigma(k,1)$ in $H_T^*(LG_n)$:
\begin{equation}
\sigma(k,1)=\sigma(k)\sigma(1)-2\sigma(k+1)-2x_{n-k+1}\sigma(k)\quad(2\leq k\leq n),
\label{tworows}
\end{equation}
where $\sigma(j)=0$ for $j>n.$
The above 
expression is a consequence of the following formula
for Ivanov's functions
\beqn
Q_{k,1}(x|a)=Q_k(x|a)Q_1(x|a)-2Q_{k+1}(x|a)
-2a_{k+1}Q_k(x|a).\label{Qk1}
\eeqn
As illustrated by this example, 
we need a correction term 
to classical formula (\ref{quad}).

\bigskip

To generalize (\ref{tworows}), we prove the next proposition,
which is also interesting from a purely combinatorial 
point of view.
In this section,
$x=(x_1,x_2,\ldots)$ and $a=(a_2,a_3,\ldots)$
are two sequences of infinite variables.
We can define $Q_\lambda(x|a)$ for any 
strict partition $\lambda.$
They are in the ring $\Z[a_2,a_3,a_4,\ldots]\otimes_\Z \Gamma,$
where $\Gamma$ denote a distinguished subring 
spanned by the Schur's $Q$-functions in 
``the ring of symmetric functions $\Lambda$'' (\cite{M1}).
For the detail of definition for $Q_\lambda(x|a)$, see \cite{Iv2} and also Section \ref{App}.
Note that, if we substitute $x_j=0\;(j>n)$ 
for $Q_\lambda(x|a)\;(\lambda\in \mathcal{SP}_n)$ we 
can recover the polynomial introduced by Definition \ref{nimmo}.
Let $h_r$ (resp. $e_r$) denote the $r$-th 
complete (resp. elementary) symmetric function.

\begin{prop}\label{expansion}
Let $k>\ell>0.$ We have
\begin{equation}
Q_{k,\ell}(x|a)
=Q_k(x|a)Q_\ell(x|a)+2\sum_{i=1}^\ell(-1)^i
Q_{k+i}(x|a)Q_{\ell-i}(x|a)
+G_{k,\ell}(x|a),\label{exp1}
\end{equation}
where 
\begin{equation}
G_{k,\ell}(x|a)=\sum_{r=k}^{k+\ell-1}\sum_{s=0}^{k+\ell-1-r}
f_{k,\ell}^{r,s}(a)
Q_{r}(x|a)Q_{s}(x|a),\label{exp2}
\end{equation}
and the coefficient $f_{k,\ell}^{r,s}(a)$ is given by
\begin{equation}
f_{k,\ell}^{r,s}(a)=(-1)^{\ell-s}
\sum_{j=0}^{k+\ell-r-s}
2\,h_{k+\ell-r-s-j}(a_{k+1},a_{k+2},\ldots,a_{r+1})\,e_{j}
(a_{s+2},\ldots,a_{\ell-1},a_{\ell}).\label{exp3}
\end{equation}
\end{prop}

{\bf Proof.} 
We use the equation (8.2) of \cite{Iv2} that reads
\begin{equation}
Q_{k+1,\ell}+Q_{k,\ell+1}+(a_{k+1}+a_{\ell+1})Q_{k,\ell}
=Q_{k}Q_{\ell+1}-Q_{k+1}Q_{\ell}
+(a_{\ell+1}-a_{k+1})Q_kQ_\ell,\label{Iv8.1}
\end{equation} for $k>\ell >0,$
where we denote $Q_{r,s}(x|a)$ simply by $Q_{r,s}.$
By this equation, it is easy to see that each function $Q_{k,\ell}$ is
a linear combination of the functions $Q_rQ_s\,(r>s\geq 0).$
Note that the functions $Q_rQ_s\,(r>s\geq 0)$ are
linearly independent over the ring $\Z[a_2,a_3,\ldots]$
(This fact can be seen from Prop. 2.11 in \cite{Iv2} and \cite{M1} III, (8.9)).

We shall prove the proposition by induction on $\ell.$
The case $\ell=1$ is true by (\ref{Qk1}). 
Let $\ell>1.$ Suppose the proposition holds for $\ell.$
We have an expansion 
\begin{equation}
Q_{k,\ell+1}=\sum_{r>s\geq 0}g_{k,\ell+1}^{r,s}(a)Q_rQ_s,
\end{equation}
with coefficients $g_{k,\ell+1}^{r,s}(a)\in \Z[a_2,a_3,\ldots].$
Our task is to show $g_{k,\ell+1}^{r,s}=f_{k,\ell+1}^{r,s}.$
By extracting the coefficient of $Q_kQ_\ell$ in 
both hand sides of (\ref{Iv8.1}), we have
$$
g_{k,\ell+1}^{k,\ell}
+(a_{\ell+1}+a_{k+1})=
a_{\ell+1}-a_{k+1}.
$$
Hence we have
$g_{k,\ell+1}^{k,\ell}=-2a_{k+1}=f_{k,\ell+1}^{k,\ell}.$
Let $(r,s)\ne (k,\ell)$ with $r+s<k+\ell.$
By comparing the coefficients of $Q_rQ_s$ 
in 
both hand sides of (\ref{Iv8.1}), we have
\begin{equation}
f_{k+1,\ell}^{r,s}+g_{k,\ell+1}^{r,s}+(a_{k+1}+a_{\ell+1})f_{k,\ell}^{r,s}=0.
\end{equation} 
We shall prove $g_{k,\ell+1}^{r,s}=f_{k,\ell+1}^{r,s}$
by showing
\begin{equation}
f_{k+1,\ell}^{r,s}+f_{k,\ell+1}^{r,s}+(a_{k+1}+a_{\ell+1})f_{k,\ell}^{r,s}=0.
\end{equation}
This follows from the following equality:
$$
\frac{\prod_{\alpha=s+2}^\ell(1+a_\alpha z)}{\prod_{\beta=k+2}^{r+1}(1-a_\beta z)}
+z(a_{k+1}+a_{\ell+1})
\frac{\prod_{\alpha=s+2}^\ell(1+a_\alpha z)}{\prod_{\beta=k+1}^{r+1}(1-a_\beta z)}
=\frac{\prod_{\alpha=s+2}^{\ell+1}(1+a_\alpha z)}{\prod_{\beta=k+1}^{r+1}(1-a_\beta z)}.
$$
$\square.$

For example we have
\begin{eqnarray*}
G_{k,1}&=&-2a_{k+1}Q_k,\\
G_{k,2}&=&
2(a_{k+1}+a_{k+2}+a_2)Q_{k+1}-2a_{k+1}Q_kQ_1+2(a_{k+1}^2+a_2a_{k+1})Q_k,\\
G_{k,3}&=&-2(a_{k+1}+a_{k+2}+a_{k+3}+a_2+a_3)Q_{k+2}
+2(a_{k+1}+a_{k+2}+a_3)Q_{k+1}Q_1\\
&&-2\left(a_{k+1}^2+a_{k+1}a_{k+2}+a_{k+2}^2+(a_{k+1}+a_{k+2})(a_2+a_3)+a_2a_3\right)Q_{k+1}
\\
&&-2a_{k+1}Q_kQ_2+2(a_{k+1}^2+a_{k+1}a_3)Q_kQ_1
-2\left(a_{k+1}^3+a_{k+1}^2(a_2+a_3)+a_{k+1}a_2a_3\right)Q_k,
\end{eqnarray*}
where we denote $Q_{r,s}(x|a)$ simply by $Q_{r,s}.$

Proposition \ref{expansion} combined with Theorem \ref{res}
give rise to a polynomial expression for $\sigma(r,s)$ with $n\geq r>s>0$
in terms of the special classes $\sigma(k)\;(1\leq k\leq n).$
For example, we have
\begin{eqnarray*}
\sigma(k,2)&=&\sigma(k)\sigma(2)-2\sigma(k+1)\sigma(1)+2\sigma(k+2)
-2x_{n-k+1}\sigma(k)\sigma(1)\\
&+&2(x_{n-k+1}+x_{n-k}+x_n)\sigma(k+1)
+2(x_{n-k+1}^2+x_{n-k+1}x_n)
\sigma(k)
\end{eqnarray*}
for $2<k\leq n$, with $\sigma(j)=0$ for $j>n$ (cf. Proposition \ref{vanish}).

\bigskip
The next proposition will be used in Section \ref{Ring}.
\begin{prop}\label{Xii}
For $k\geq 1$, we have
\beqn
{Q}_k(x|a)^2
+2\sum_{i=1}^k (-1)^i 
{Q}_{k+i}(x|a){Q}_{k-i}(x|a)+
\sum_{r=k}^{2k-1}\sum_{s=0}^{2k-1-r}
{f_{k,k}^{r,s}}(a){Q}_r(x|a){Q}_s(x|a)=0.\label{QkSq}
\eeqn
\end{prop}
{\bf Proof.}
The proof Lemma \ref{expansion}  
is valid also 
for $k=\ell$ with $Q_{k,k}(x|a)=0$ for $k\geq 1.$
$\square.$

\section{Presentation of the ring $H_T^*(LG_n)$}\label{Ring}
\setcounter{equation}{0}
As an application of Theorems \ref{res} and \ref{Giam},
we obtain a presentation of the ring $H_T^*(LG_n)$ 
in terms of generators and relations.
Consider the ring 
$\Z[a]=\Z[a_2,a_3,\ldots,a_{n+1}].$
Throughout the section, 
we identify $\mathcal{S}=\Z[x_1,\ldots,x_n]$ and $\Z[a]$ 
by the isomorphism 
$\iota_n:\Z[a]\rightarrow\mathcal{S}$ of rings
given by 
\begin{equation}
\iota_n(a_{j})=x_{n-j+2}\,(2\leq i\leq n).\label{isom}
\end{equation}

\subsection{Statement of the result}

Let $X_1,\ldots,X_n$ denote a set of indeterminates.
Set $X_0=1$ and $X_j=0$ for $j>n$ (cf. Proposition \ref{vanish}).
Let $k,\ell$ be $n\geq k\geq \ell\geq 0.$
Consider the following 
elements of the polynomial ring $\mathcal{S}[X_1,\ldots,X_n]:$
\begin{equation}
X_{k,\ell}=X_kX_\ell
+2\sum_{i=1}^{\mathrm{min}(n-k,\ell)} (-1)^i 
X_{k+i}X_{\ell-i}+
\sum_{r=k}^{\mathrm{min}(n,k+\ell-1)}\sum_{s=0}^{k+\ell-1-r}
{f_{k,\ell}^{r,s}}(a)\,X_rX_s,\label{Xkl}
\end{equation}
where $f_{k,\ell}^{r,s}(a)$ is given 
by the right hand side of (\ref{exp3}).
Since we restrict $r\leq n$,
we can consider $f_{k,\ell}^{r,s}(a)$ 
to be in $\mathcal{S}$ via the isomorphism $\iota_n.$
Note also that we also consider 
the case of $\ell=k$.
Define an ideal $\mathcal{I}_n=\langle X_{1,1},\ldots,X_{n,n}\rangle$ and
consider the quotient 
ring $$\mathcal{R}_n=\mathcal{S}[X_1,\ldots,X_n]/\mathcal{I}_n.$$

We shall define a morphism of $\mathcal{S}$-algebras
$\phi:\mathcal{R}_n
\longrightarrow 
H_{T}^*(LG_n)$ by setting $X_i$ to $\sigma(i)\,(1\leq i\leq n).$ 
\begin{lem}
The map $\phi$ is well-defined.
\end{lem}
{\bf Proof.}
Define a morphism of $\mathcal{S}$-algebras 
$\tilde{\phi}:\mathcal{S}[X_1,\ldots,X_n]\rightarrow H_T^*(LG_n)$ by
$\tilde{\phi}(X_i)=\sigma(i).$
For $k$ with $1\leq k\leq n$, and $\mu\in W^P$, we have
\beqn
&&\tilde{\phi}(X_{k,k})|_{\mu}\nonumber\\
&=&\sigma(k)^2|_{\mu}
+2\sum_{i=1}^{\mathrm{min}(n-k,k)} (-1)^i 
\sigma({k+i})|_{\mu}\sigma({k-i})|_{\mu}+
\sum_{r=k}^{\mathrm{min}(n,2k-1)}\sum_{s=0}^{2k-1-r}
\iota_n({f_{k,k}^{r,s}}(a))\,\sigma(r)|_{\mu}\sigma(s)|_{\mu}\nonumber\\
&=&Q_k(x_\mu|x_{\langle n\rangle})^2
+2\sum_{i=1}^{k} (-1)^i 
Q_{k+i}(x_\mu|x_{\langle n\rangle})Q_{k-i}(x_{\mu}|x_{\langle n\rangle})\nonumber\\
&+&
\sum_{r=k}^{2k-1}\sum_{s=0}^{2k-1-r}
{f_{k,k}^{r,s}}(x_{\langle n\rangle})\,Q_r(x_\mu|x_{\langle n\rangle})
Q_s(x_\mu|x_{\langle n\rangle})\nonumber
\eeqn
where in the second equality,
we used Theorem \ref{res} and
a vanishing property (Proposition \ref{vanish}).
We can see the last expression is zero
by specializing (\ref{QkSq}) (see 
Definition \ref{specialization}).
Thus we have $\tilde{\phi}(X_{k,k})|_\mu=0$ for all $\mu\in W^P$,
and hence $\tilde{\phi}(X_{k,k})=0.$
So $\tilde{\phi}$ induces $\phi:\mathcal{R}_n\rightarrow H_T^*(LG_n)$ such
that $\phi(X_i)=\sigma(i)\;(1\leq i\leq n).$
$\square.$

\bigskip

\begin{Def}
Let $\lambda=(\lambda_1>\cdots>\lambda_{2r}\geq 0)$ be in $\mathcal{SP}_n.$
We introduce the following Schur-type Pfaffian 
$$
X_\lambda=\mathrm{Pf}(X_{\lambda_i,\lambda_j})_{1\leq i,j\leq 2r}.
$$
\end{Def}

\begin{thm}\label{RingPresen}
 There exists an
isomorphism of $\mathcal{S}$-algebras:
$$\phi:
\mathcal{R}_n\longrightarrow 
H_T^*(LG_n)
$$
sending $X_i$ to $\sigma(i)\;(1\leq i\leq n)$
and the Pfaffian $X_\lambda$ to the equivariant Schubert class $\sigma(\lambda).$
\end{thm}

\bigskip
By definition of $\phi$
and Giambelli formula (\ref{Giambelli}), we have $\phi(X_\lambda)=\sigma(\lambda).$
Moreover, since $\sigma(\lambda)\;(\lambda\in \mathcal{SP}_n)$ 
generates $H_T^*(LG_n)$ as an $\mathcal{S}$-module,
$\phi$ is surjective.
The rest of this section is devoted to the 
proof 
of injectivity of $\phi.$

\subsection{A monomial ordering}

Here we give a preliminary discussion to 
prove Theorem \ref{RingPresen}. The argument below is 
quite similar to the one in Macdonald \cite{M1} (III, 8),
however a different ordering 
on the partitions will be used,
which proves to be
useful in our situation.

For any partition $\lambda=(1^{e_1}2^{e_2}\cdots n^{e_n})$,
we set 
$$
X^\lambda=X_{1}^{e_1}\cdots X_n^{e_n}.
$$
By $\mathrm{deg}(\lambda)$ we denote the
degree $\sum_{i=1}^ne_i$ of the monomial
$X^\lambda.$
Let $\mu=(1^{e_1'}\cdots n^{e_n'})$ be another partition.
We write $\lambda\succ\mu$ if
$\deg(\lambda)>\deg(\mu),$ or 
$$\deg(\lambda)=\deg(\mu)\;
\mbox{and there is}\; k \;\mbox{such that}\;
e_1=e_1',\ldots,e_{k}=e_{k}'\;\mbox{and} 
\;e_{k+1}<e_{k+1}'.$$
Then we also write $X^\lambda\succ X^\mu.$
This is a monomial ordering 
called the {\it grevlex order}
with
$X_1\prec X_2\prec\cdots\prec X_n.$ In particular,
if we have $\lambda\succ\mu,$ then $\lambda+\nu\succ \mu+\nu$ for
any partition $\nu.$

\begin{lem}\label{q_gen}
Let $\lambda=(1^{e_1}2^{e_2}\cdots n^{e_n})$ be a
partition.
If $\lambda$ is not strict, 
then $X^\lambda$ is an
$\mathcal{S}$-linear combination
of the $X^\mu$ with $\mu\in \mathcal{SP}_n,$
and $\mu\prec\lambda.$
In particular, the monomials
$X^\lambda\;(\lambda\in \mathcal{SP}_n)$ generate
$\mathcal{R}_n$ as an $\mathcal{S}$-module.
\end{lem}

{\bf Proof.} First note that if $\lambda$ is strict then we have $\lambda\in \mathcal{SP}_n.$
We prove the first statement by induction, assuming the claim for all 
partition $\mu$ such that $\mu\prec\lambda.$
If $\lambda$ is not strict then for some $k$ we have 
$e_k\geq 2.$
We have the following relation:
\begin{equation}
X_k^2
=-2\sum_{i=1}^{\mathrm{min}(n-k,k)} (-1)^i 
X_{k+i}X_{k-i}-
\sum_{r=k}^{\mathrm{min}(n,2k-1)}\sum_{s=0}^{2k-1-r}
{f_{k,k}^{r,s}}(a)X_rX_s.\label{Xsq}
\end{equation}
We can see that the monomials appearing in the 
right hand side of the above equation is strictly lower than $X_k^{2}$
in the grevlex order $\prec$.
Replacing the factor $X_k^2$ in $X^\lambda$ by the right hand side
of (\ref{Xsq}), we can express $X^\lambda$ as an $\mathcal{S}$-linear
combination of the $X^\mu$'s where each $\mu$ is a
partition such that $\mu\prec \lambda.$
By the inductive hypothesis the claim is true for each $X^\mu,$
and the proof completes.
$\square.$

\begin{lem}\label{Qq} Let $\lambda\in \mathcal{SP}_n.$
The Pfaffian $X_\lambda$ is written in the form
$$
X_\lambda=X^\lambda+\sum_{\mu}b_{\lambda\mu}(a)X^\mu
$$
with coefficients $b_{\lambda\mu}(a)\in \mathcal{S},$
where the sum is over $\mu\in \mathcal{SP}_n$ such that $\mu\prec\lambda.$
\end{lem}

{\bf Proof.} Let $\lambda=(\lambda_1>\cdots>\lambda_{2r}\geq 0)$ be
a strict partition in $\mathcal{SP}_n.$
We proceed by induction on $r.$
Let $r=1.$
If $\lambda=(i)$ with $1\leq i\leq n$ the Lemma is clear.
For two-row type the Lemma is true by (\ref{Xkl}).
Let $r\geq 2$ and assume
the Lemma holds for all $\mu=(\mu_1>\cdots>\mu_{2s}\geq 0)
\in \mathcal{SP}_n$ with $s<r.$
From the definition of the Pfaffian
it follows that
$$
X_{\lambda}=\sum_{j=2}^{2r}(-1)^{j}X_{\lambda_1,\lambda_j}
X_{\lambda_2,\cdots,\widehat{\lambda_j},\cdots,\lambda_{2r}}.
$$
By the inductive hypothesis, 
we have
$$
X_{\lambda_2,\cdots,\widehat{\lambda_j},\cdots,\lambda_{2r}}=X_{\lambda_2}\cdots\widehat
{X_{\lambda_j}}\cdots X_{\lambda_{2r}}+F_j
$$
where $F_j$ is a $\mathcal{S}$-linear combination
of $X^\mu$'s with $\mu\in \mathcal{SP}_n$ such 
that $\mu\prec(\lambda_2,\cdots,\widehat{\lambda_j},\cdots,\lambda_{2r}).$
Then it is easy to see that the Lemma holds for $\lambda.$
$\square.$

\bigskip

From Lemmas \ref{q_gen} and \ref{Qq}, we have
the following.
\begin{lem}\label{Q_gen}
The Pfaffians $X_\lambda\;(\lambda\in \mathcal{SP}_n)$ generate
$\mathcal{R}_n$ as an $\mathcal{S}$-module.
\end{lem}

\subsection{ Completion of the proof of Theorem \ref{RingPresen}.} 
It remains to 
prove the injectivity of $\phi.$
 Let $F$ be in $\mathrm{Ker}(\phi).$
By Lemma \ref{Q_gen} we have
$$
F=\sum_{\lambda\in \mathcal{SP}_n}c_\lambda(a)X_\lambda
$$
with coefficients $c_\lambda(a)\in \mathcal{S}.$
We know $\phi(X_\lambda)=\sigma(\lambda).$
So we have 
$0=\sum_\lambda \iota_n(c_\lambda(a))\sigma(\lambda).$ 
Since $\sigma(\lambda)$ are linearly independent over $\mathcal{S}$,
$\iota_n(c_\lambda(a))=0$ for all $\lambda\in \mathcal{SP}_n.$
Hence we have $c_\lambda(a)=0\,(\lambda\in \mathcal{SP}_n)$ and
$F=0.$
\bigskip

\section{Appendix}\label{App}
\setcounter{equation}{0}

For the reader's convenience, we provides
a summary of some properties of 
$Q_\lambda(x|a).$
We also prove a vanishing property (Proposition \ref{vanish})
essentially used in the main body of the paper.

We use standard notation for symmetric
functions as in Macdonald's book \cite{M1}.
Let $\Lambda$ denote 
the {\it ring of symmetric functions} in 
infinitely many indeterminates $x=(x_1,x_2,\ldots).$
The ring $\Lambda$ is graded as $\Lambda=\oplus_{k= 0}^\infty \Lambda^k$
and each graded part $\Lambda^k$ has a $\Z$-basis
consisting of the
{\it monomial symmetric functions} $m_\lambda=m_\lambda(x)$
(for all partitions $\lambda$ of $k$).

Recall an expression for the $Q_k(x)$
the Schur's $Q$-functions
for the one-row partition
$$
Q_k(x)=\sum_{\lambda}2^{\ell(\lambda)}m_\lambda(x),
$$
where the sum runs over the all partitions $\lambda$ of $k$ and 
$\ell(\lambda)$ is the length of 
$\lambda$, the number of nonzero parts of $\lambda.$
Let $\Gamma$ be the subring of $\Lambda$ generated by $Q_k$:
$$
\Gamma=\Z[Q_1,Q_2,Q_3,\ldots].
$$
We have a gradation $\Gamma=\oplus_{k=0}^\infty\Gamma^k$ where 
$\Gamma^k=\Gamma\cap \Lambda^k.$
The Schur's $Q$-functions $Q_\lambda(x)$, with $\lambda$ srict
partition of $k$, form a distinguished $\Z$-basis of $\Gamma^k.$ 

Let $a_2,a_3,a_4,\ldots$ be an infinite sequence of 
independent variables. We set $a_1=0.$
Ivanov introduced a factorial analogue of $Q$-functions $Q_\lambda(x|a)$
defined for any strict partition $\lambda.$ 
Each $Q_\lambda(x|a)$ is an element of the ring $\Z[a_2,a_3,\ldots]\otimes_\Z\Gamma.$
In particular, we have, by Ivanov \cite{Iv2}, Theorem 8.2, 
$$
Q_k(x|a)=\sum_{j=0}^{k-1}(-1)^je_{j}(a_2,a_3,\ldots,a_k)Q_{k-j}(x).
$$
For $k>\ell>0,$ we can {\it define} $Q_{k,\ell}(x|a)$ by Proposition \ref{expansion}.
Moreover, for arbitrary strict 
partition $\lambda,$ we
have
\beqn
Q_\lambda(x|a)=\mathrm{Pf}(Q_{\lambda_i,\lambda_j}(x|a))_{1\leq i<j\leq 2r},
\label{IvSchur}
\eeqn
where we write
$\lambda=(\lambda_1,\lambda_2,\ldots,\lambda_{2r})$ with
$\lambda_1>\cdots>\lambda_{2r}\geq 0.$

\bigskip

The following result is 
very important. Recall Definition \ref{specialization}
for the meaning of $Q_\lambda(x_\mu|x_{\langle n\rangle}).$

\begin{prop} \label{vanish}
$Q_\lambda(x_\mu|x_{\langle n\rangle})$ 
vanishes identically unless $\lambda\in \mathcal{SP}_n.$
\end{prop}
{\bf Proof.}
We prove the Proposition for $\lambda=(k).$
In \cite{Iv2}, Ivanov derived the following 
equation (Theorem 8.2 in \cite{Iv2})
\begin{equation}
\sum_{k=0}^\infty
\frac{Q_k(x|a)z^k}{\prod_{j=1}^k(1-a_{j+1}z)}
=\prod_{i=1}^\infty \frac{1+x_iz}{1-x_iz}.\label{onerow}
\end{equation}
If we specialize the variables as in the statement of the Proposition,
we have
$$
\sum_{k=0}^n
\frac{Q_k(x_\mu|x_{\langle n\rangle})z^k}{\prod_{j=1}^k(1-x_{n+1-j}z)}
+
\frac{\sum_{k>n}Q_k(x_\mu|x_{\langle n\rangle})z^k}{\prod_{j=1}^n(1-x_{j}z)}
=\prod_{1\leq i\leq n,\delta_i=1} \frac{1+x_iz}{1-x_iz},
$$
where $(\delta_1,\ldots,\delta_n)\in \mathcal{M}_n$ corresponds
to $\mu\in W^P.$
Multiplying $\prod_{1\leq i\leq n}(1-x_iz)$
to both hand sides,
we have
$${\sum_{k>n}Q_k(x_\mu|x_{\langle n\rangle})z^k}
=-\sum_{k=0}^n
{Q_k(x_\mu|x_{\langle n\rangle})z^k}{\prod_{j=1}^{n-k}(1-x_{j}z)}
+\prod_{i=1}^n {(1+(-1)^{\delta_i+1}x_iz)}.
$$
The right hand side of the equation is a polynomial in $z$ of
degree lower than $n$ and we are done.
For general $\lambda,$ the Proposition follows from
 Proposition \ref{expansion}
and the Pfaffian formula (\ref{IvSchur}) for $Q_\lambda(x|a)$.
$\square.$

\begin{small}
{\scshape Department of Applied Mathematics,
Okayama University of Science,
Okayama 700-0005, JAPAN}
\end{small}

\textit{E-mail address}: \tt{ike@xmath.ous.ac.jp}

\end{document}